\documentclass[epsbox,12pt,reqno]{amsart}
\usepackage{amsmath,amssymb,latexsym,amscd,amsthm} % ams auxiliaries.
\usepackage{graphicx}

\theoremstyle{definition}

\theoremstyle{remark}

\numberwithin{equation}{section}

\title{The Colourful Feasibility Problem}

\author{Antoine Deza}
\author{Sui Huang}
\address{Advanced Optimization Laboratory,
Department of Computing and Software,
McMaster University,
Hamilton, Ontario,
Canada
L8S 4K1.
}
\email{\{deza,huangs3,terlaky\}@mcmaster.ca}
\author{Tamon Stephen} 
\address{
Department of Mathematics,
Simon Fraser University,
8888 University Drive,
Burnaby, British Columbia, Canada  V5A 1S6.
}
\email{tamon\_stephen@sfu.ca}
\author{Tam{\' a}s Terlaky}

\dedicatory{Dedicated to Leonid Khachiyan.}


\begin{document}

%%%%% Insert abstract here.

\begin{abstract}
We study a colourful generalization of the linear programming
feasibility problem, comparing the algorithms introduced by 
B{\' a}r{\' a}ny and Onn with new methods.  
We perform benchmarking on generic and ill-conditioned problems,
as well as as recently introduced highly structured problems.  
We show that some algorithms can lead to cycling or slow convergence,
but we provide extensive numerical experiments which show that
others perform much better than 
predicted by complexity arguments.  We conclude that the most
efficient method for all but the most ill-conditioned problems
is a proposed multi-update algorithm.
\end{abstract}

\newif\ifapps        % Include appendices 3 and 4?
\appstrue

%%%%%%%%%%%%%%%%%%%%%%%%%%%%%%%%%%%%%% Insert body here.

%%%%%%%%%%%%% November 2005 %%%%%%%%%%%
%
%% Include file for main template.
%

%%% Special commands used.
\def\sd{\operatorname{depth}}     % Simplicial depth.
\def\csd{\mathbf{depth}}          % Colourful simplicial depth.
\def\S{\mathbf{S}}                % Configuration.
\def\R{\mathbb{R}}
\def\core{\operatorname{core}}
\def\cone{\operatorname{cone}}
\def\conv{\operatorname{conv}}
\def\new{\operatorname{new}}
\def\aff{\operatorname{Aff}}
\def\interior{\operatorname{int}}
\def\bara{B\'ar\'any}
\def\baraonn{B\'ar\'any and Onn}
\def\cara{Carath\'eodory}
\def\Sph{\mathbb{S}}              % Sphere.
\def\zero{{\bf 0}}
\def\Deltab{\mathbf{\Delta}}
\def\B{{\bf A1}}                       % Algorithms
\def\BO{{\bf A2}}
\def\Bp{{\bf A3}}
\def\BOp{{\bf A4}}
\def\H{{\bf A5}}
\def\MV{{\bf A6}}
\def\Rand{{\bf A7}}
\def\A{{\bf Ax}}
\def\rg{{\bf G1}}                      % Generators
\def\tz{{\bf G2}}
\def\to{{\bf G3}}
\def\Sp{{\bf G4}}
\def\Sd{{\bf G5}}
\def\Sm{{\bf G6}}
\def\Sdelta{$\mathbf{S^{\Delta}}$}
\def\Splus{$\mathbf{S^{+}}$}
\def\Sminus{$\mathbf{S^{-}}$}

\maketitle

%
% INTRODUCTION
%
\section{Introduction}
Given colourful sets $S_1, \ldots, S_{d+1}$ of points in $\R^d$
and a point $p$ in $\R^{d}$, the {\it colourful feasibility problem} 
is to express $p$ as a convex combination of
points $x_1, \ldots, x_{d+1}$ with $x_i \in S_i$ for each $i$.
This problem was presented by \bara\ in 1982 \cite{Bar82}.
The monochrome version of this problem, expressing $p$ as a 
linear combination of points in a set $S$, is a traditional 
linear programming feasibility problem.

In this paper, we study algorithms for the colourful feasibility problem
with a core condition from an experimental point of view.  
We learn several things.
First this problem is easy in a practical sense --
we expend more effort to generate difficult examples than to solve
them.  
Second, while the classical algorithms for this problem already perform
quite well, we introduce modifications that achieve a substantial improvement
in practical performance.
Third, we construct examples where ill-conditioning leads
to slow convergence for the some otherwise very effective
algorithms.
And finally, we remark that a simple greedy heuristic 
provides competitive results in practice  but we find a case where
it fails to solve the problem at all.
Additionally we provide benchmarking that we hope will encourage
research on this attractive problem.

\section{Definitions and Background}
We concentrate on the important subcase of colourful feasibility problem
where we have $d+1$ points of
each colour, and $p \in \conv(S_i)$ for $i=1,\ldots,d+1$.  
We call $\bigcap_{i=1}^{d+1} \conv(S_i)$ the {\it core} of the configuration.
We will call such a problem a {\it colourful feasibility} problem,
in this paper colourful feasibility problems are 
assumed to have a non-empty core.
In this core case, by \bara's colourful \cara\ theorem \cite{Bar82}, 
a solution is guaranteed to exist, and the problem is to exhibit 
a solution.  
Recently \bara's result has been strengthened to show that quadratically
many solutions must exist, see \cite{BM06} and \cite{ST06}.
The problem of finding a solution to a colourful feasibility problem
is described
in \cite{BO97} as ``an outstanding problem on the border line
between tractable and intractable problems''.

Several close relatives of the colourful feasibility problem
are known to be difficult.
For example, the case where we have $d$ colours in $\R^d$ and no restriction
on the size of the sets has been shown to be strongly NP-complete
through a reduction of 3-SAT.  We refer to \cite{BO97} for more details.

In \cite{Bar82}, \bara\ proposed a finite algorithm \B\ to solve 
colourful feasibility,
and in \cite{BO97} \baraonn\ analyzed the complexity of \B\ 
and a second algorithm \BO.  
Both these algorithms are essentially geometric, and 
the complexity guarantees depend crucially on having the 
point $p$ in the {\it interior} of the core.  
In effect, the distance between $p$ and the boundary 
of the core can be considered as a measure of the conditioning 
of the problem.  Thus for a configuration
$\S$ we define $\rho$ to be the radius of the largest ball around 
$p$ that is contained in the core. % $\bigcap_{i=1}^{d+1} \conv(S_i)$.
The results for \B\ and \BO\ are effectively that they are
polynomial in $d$ and $1/\rho$.
We remark, though, that for configurations of $d+1$
points in $d+1$ colours on the unit sphere $\Sph^d \subseteq \R^d$, 
$\rho$ will be small even if
the problem has a favourable special structure, and quite small
otherwise. 

Without loss of generality, we can take the point $p$ to
be the vector $\zero$ in $\R^d$.
Some additional preprocessing will be helpful.
If $\zero$ is a point
in one of the $S_i$'s, then the solution to the colourful 
feasibility problem is trivial.  Otherwise, we can % without loss of generality 
scale the points of the $S_i$'s
so that they lie on the unit sphere $\Sph^d$.  The coordinates in
any resulting convex combination can then be unscaled as a 
post-processing step.

We call a system of $d+1$ sets of $d+1$ points a {\it configuration},
and often denote it as $\S=\{S_1,\ldots,S_{d+1}\}$.  
We use a the bold font to signal a colourful object, except with
$\zero$ where bold is used to distinguish the vector from a scalar.
We remark that restricting the sets to have size $d+1$
is not a burden since, given a larger set, solving a monochrome linear 
feasibility problem allows us to efficiently find a basis of size $d+1$ 
with $\zero$ in its convex hull.

%
% ALGORITHMS
%
\section{Seven Algorithms}\label{se:algs}
In this paper we consider the theoretical and practical performance
of seven algorithms for finding a colourful basis.
The algorithms considered are the algorithms of \bara\ 
\B\ and \baraonn\ \BO, modifications of these algorithms which
update multiple colours at each stage, which we will call \Bp\ and \BOp\,
and a hybrid \H\ of these designed to take advantage of the strengths of
both algorithms.  For purposes of comparison, we also consider
two simple approaches that perform well under certain 
circumstances: 
a greedy heuristic where we choose the adjacent simplex of
maximum volume \MV\  
and a random sampling approach \Rand.
All our implementations are initialized with using the first points from 
each colour.
% Sui Huang has written MATLAB code to implement each of these algorithms.
% We have made this freely available on-line \cite{Hua05}.
Following are descriptions of the algorithms,
see \cite{Hua05} for MATLAB implementations of each.
Besides \Rand, they are implemented as
pivoting algorithms with the respective pivot selection rule. 

\subsection{\bara's Algorithm \B}\label{se:b}
% We now offer brief descriptions of the algorithms.
We begin with the algorithm proposed by \bara\ \cite{Bar82},
which is a pivoting algorithm.  It begins with say a random colourful
simplex $\Deltab$.  The point $x$ nearest to $\zero$ in $\Deltab$ is
computed.  If $x \ne \zero$, then $x$ must lie on some facet of $\Deltab$.  
Consider the colour $i$ of the vertex of $\Deltab$ that is not
on this facet.  Look for
the point $t$ of colour $i$ minimizing the inner product 
$\langle t,x \rangle$.  Then we replace the point of colour $i$
from $\Deltab$ with the point $t$ to get a new simplex.
The algorithm then repeats beginning with the new simplex.

The convergence of this algorithm relies on the fact that $\zero$ 
is in the core of the
configuration.  For this reason the affine hyperplane perpendicular
to the vector $x$ cannot
separate $\zero$ from the points of colour $i$.  Thus the next simplex
will have a point closer to $\zero$ than $\Deltab$ did, and the algorithm
will converge in finitely many steps.  If, additionally, the core
has radius at least $\rho$ around $\zero$, then there is a guarantee 
on the amount of progress in a given step, which depends on $\rho$.
Effectively the guarantee is that the number of iterations of \B\ 
is $O(1/\rho^2)$.  Since an iteration can be done in polynomial
time, this proves that \B\ runs in time polynomial in the input data
and $1/\rho$.  
Consult \cite{BO97} for details and a proof.  

We note that the complexity of a single iteration is dominated
by the cost of the nearest point subroutine.  
This is can be solved as a continuous optimization problem,
but complicates our life with numerical issues: It can be
solved to less or greater precision, either risking numerical
error or increasing the running time.  For the purposes of
our benchmarking, we used the {\tt MATLAB} built-in {\tt quadprog()}
which gave fairly good results, 
see Section~\ref{se:tpi}. 

\subsection{\baraonn's Algorithm \BO}\label{se:bo}
The reliance of \B\ on nearest point calculations is certainly a
disadvantage.
Partly motivated by this, \baraonn\ proposed an alternate
algorithm for the colourful feasibility problem whose calculations involve 
only linear algebra.
This algorithm, \BO, is described in \cite{BO97}.

Essentially, the closest point $x$ to $\zero$ on 
the simplex $\Deltab$ is replaced in this algorithm by 
a point $y$ on the boundary of $\Deltab$ that can be computed 
algebraically. % as a convex combination of the vertices of $\Deltab$.  
The initial choice of $y$ could be one of the vertices of the initial
simplex.  In subsequent iterations,
a colour $j$ corresponding to a zero coefficient in $y$ is chosen.
An improving vertex $v$ of colour $j$ is found, and $y_{\new}$ 
is updated by projecting $\zero$ onto the line segment between
$y$ and $v$ and finding where the resulting vector enters the
new simplex.  As with \B, this algorithm takes $O(1/\rho^2)$
iterations, and hence is polynomial in the input data and $1/\rho$,
see \cite{BO97}.

The implementation of \BO\ proposed in \cite{BO97} takes time
$\Theta(d^4)$ for a single iteration.  The bottleneck is
computing $y_{\new}$, which is the intersection of the line segment 
from $\zero$ to 
a point $p$ and the new simplex.  In fact we observe that this 
can be done in time $O(d^3)$. 
First, compute the defining equations for the simplex $Ay_{\new} \ge b$
by inverting the homogenized matrix of the vertices.
We know the intersection point will be of the form $y_{\new}= \alpha p$. 
We can substitute this into the above inequalities to get 
$\alpha (Ap) \ge b$ and simply take $\alpha$ to be the maximum
value of $b_i/A_i p$ for $i=1,2, \ldots, d+1$.  
This is implemented in \cite{Hua05}.  

\subsection{Multi-update \bara\ \Bp}\label{se:bp}
We are interested in getting practically effective algorithms
for the colourful feasibility problem.  
To that end, we propose the following modification of \B.  
If it happens that the nearest point $x$ to $\zero$ of the current
simplex $\Deltab$ lies on a lower-dimensional face of $\Deltab$
- i.e.,~on more than one facet -
then we update {\it every} colour that is not a vertex of that face
before recomputing $x$.  Since all the new points will be on the
$\zero$ side of hyperplanes separating $\zero$ and $\Deltab$
through $x$, the convergence proofs of \B\ and \BO\ still apply
to this algorithm.  The advantage of this new algorithm, which
we call \Bp, is that when possible it updates several colours
without recomputing a nearest point.

Since this algorithm makes at least as much progress as \B\ at
each iteration, we get convergence in at most the same
number of iterations.  A given iteration may take longer,
since it has to update multiple points.  However, aside
from the nearest point calculation, all steps in an iteration
of \B\ can be performed in $O(d^2)$ arithmetic operations.  
Hence the additional
work per iteration of \Bp\ is $O(d^3)$, and the bottleneck
remains the single nearest point calculation.

\subsection{Multi-update \baraonn\ \BOp}\label{se:bop}
Similarly, we can adjust algorithm \BO\ to update $y$ only after
pivoting multiple colours in the case where $y$ lies on a 
low-dimensional face.  This is particularly useful at the
start if we use the setup proposed in \cite{BO97} where the
initial point $y$ is a vertex of $\Deltab$.
We call this algorithm \BOp.

As with \Bp, we expect this algorithm to take no more iterations
than the algorithm on which it is based, namely \BO.
Again we note that all steps in an iteration of \BO\ except
for computing the intersection of a line segment and a point 
take $O(d^2)$ arithmetic operations, 
so the additional work per iteration of \BOp\ 
as compared to \BO\ is at most $O(d^3)$.
Thus an iteration of \BOp\ will be asymptotically at most a 
constant factor slower than an iteration of \BO.

\subsection{Multi-update Hybrid \H}\label{se:h}
In Section~\ref{se:results} we describe a situation where
\BO\ and \BOp\ make extremely slow progress because they 
repeatedly return to the same simplex, see the example in
Section~\ref{se:mve}.  
A practical solution to this is to run \BOp, but
use a computationally heavy step from \Bp\ if we detect that \BOp\ 
is returning to the same simplex.  We implemented such a hybrid 
algorithm \H.

\subsection{Maximum Volume \MV}\label{se:mv}
For purposes of comparison, we also consider the performance of 
a greedy heuristic,
where we move from $\Deltab$ to an adjacent simplex of maximum 
volume given that the pivoting hyperplane separates $\Deltab$ 
from $\zero$.  This heuristic, which we call \MV, uses simpler
linear algebra than \BO, and by taking large simplices often
gets to $\zero$ in a small number of steps.

For a given candidate pivoting facet it is possible to choose
the point that generates the maximum volume simplex with that
facet by looking at the distances of the points of the candidate
colour to the hyperplane containing the facet.  
A single volume computation via a determinant can be done in
time $O(d^3)$ per candidate colour, thus an iteration of \MV\ 
takes $O(d^4)$ time.  Since the list of candidate colours may
not be all that large in typical situations, we can hope that
the cost of an iteration will often be less than that.

\subsection{Random Sampling \Rand}\label{se:rand}
Finally, we consider a very simple
{\it guess and check} algorithm where we sample simplices at
random and check to see if they contain $\zero$.  Intuitively
we would not expect such an algorithm to work well.  However, as 
discussed in \cite{DHST06}, solutions to a given colourful
feasibility problem may not be
all that rare, and in some cases can be quite frequent.
Since guessing and checking are relatively fast operations, it
worth considering the possibility that this naive algorithm 
is faster than more sophisticated algorithms at least in low
dimension.  We call this algorithm \Rand.

One attractive feature of \Rand\ is that the cost of an iteration
is low -- we only have to generate a random simplex and then
test if it contains $\zero$.  The test can be done in $O(d^3)$ time
by linear programming.

%
% EXAMPLES
%
\section{Random, Ill-conditioned and Extremal Problems}\label{se:ex}
To better understand how various algorithms perform in practice,
we produced a test suite of challenging colourful
feasibility problems, which includes generic, ill-conditioned
and highly structured problems.
In this section we describe three types of colourful feasibility
problems that we consider when evaluating the practical performance of an
algorithm.
See \cite{Hua05} for a MATLAB implementation of each of these
problem generators.

\subsection{Unstructured Random Problems}\label{se:rg}
The first class of problems we consider are unstructured random problems.
We take $d+1$ points in each of $d+1$ colours on $\Sph^d$.  The only
restriction we require is that $\zero$ is in the core
We achieve this by taking the last point to be a random
convex combination of the antipodes on $\Sph^d$ of the first $d$
points.  We call this generator \rg.

\subsection{Ill-conditioned Random Problems}\label{se:tg}
Next, we consider ill-conditioned problems.  
We place $d$ points of a given colour on the spherical cap
around the point $(0,0,\ldots,0,1)$ and the final point of that
colour in the opposite spherical cap, again as a convex combination
of the antipodes.  In our implementation of this, the maximum
angle between a chosen vector and the final coordinate axis is
a parameter, and points are concentrated towards the centre
rather than uniformly distributed on the cap.
Since the points all lie in a tube around the final coordinate
axis, we call these {\it tube} generators.  We implemented two
tube generators: \tz\ randomly places either 1 or $d$ points
of colour $i$ on the positive side of the axis, while \to\ 
always places $d$ points of colour $i$ on the positive side
of the axis.

\subsection{Problems with a Restricted Number of Solutions}\label{se:sg}
Finally, we consider problems where we control the number of
colourful simplices containing $\zero$.  The paper \cite{DHST06}
provides new bounds for the number of possible solutions to a colourful
linear program with
$\zero$ in the interior of the core.  It turns out that the
number of simplices containing $\zero$ in dimension $d$ can be
as low as quadratic in $d$,
but not lower, see \cite{BM06} and \cite{ST06}, 
or as high as $d^{d+1}+1$ (with $\rho>0$), which is more than
one third of the total number of simplices.
Constructions are given for colourful feasibility problems
attaining both these values.

The probability that a simplex generated by $d+1$ points
chosen randomly on $\Sph^d$ contains $\zero$ is $1/2^d$, 
see for example \cite{WW01}.  Thus in a uniformly generated 
random problem of the type generated by \rg, we would expect 
about $1/2^d$ of the $(d+1)^{d+1}$ colourful simplices 
to contain $\zero$.  This is not a large fraction, but in the
context of an effective pivoting algorithm such as \B\ 
which may pivot several neighbours to a given solution, and
pivot several neighbours of the first neighbour onto it, etc.,
we can entertain the idea that for a random 
configuration most simplices are close to a solution.
See Section~\ref{se:effective} for further discussion.

In any case, we would not be surprised if the difficulty
of a colourful feasibility problem increases as the number of solutions,
i.e.~simplices containing $\zero$, decreases.
To that end, we have written three problem generators
based on the constructions in \cite{DHST06}.
The first, \Sp\, generates
perturbed versions of the configuration from \cite{DHST06} 
with many solutions.
These problems have $d^{d+1}+1$ of the $(d+1)^{d+1}$ simplices
containing $\zero$, many more than random configurations, 
and we would expect them to be quite easy.
The second, \Sd, generates configurations where one point of 
each colour is close to each vertex of a regular simplex on $\Sph^d$.
There are $d!$ solutions corresponding to picking a
different colour from each vertex, note that this is
still much less than the $(d+1)^{d+1}/2^d$
expected in a random configuration.  
Finally, we have \Sm, which generates perturbed versions of the
configuration from \cite{DHST06} which has only
$d^2+1$ solutions.% in dimension $d$.  
The generators \Sp, \Sd\ and \Sm\ randomly permute the
order the points appear within each colour.

All these problems are ill-conditioned in the sense that
points are clustered closely together.  Also $\rho$ will
be quite small for \Sp\ and \Sm, although the construction
\Sd\ effectively maximizes $\rho$ for configurations on $\Sph^d$
at $1/d$.

%
% RESULTS
%
\section{Benchmarking and Results}\label{se:results}
In this section, we describe the results of computational
experiments in which we run our colourful feasibility 
algorithms against our problem generators.  
We focus on the number of iterations that an algorithm
takes to find a solution, but in Section~\ref{se:tpi}
we also include information about the cost of iterations.
The two particularly difficult, but fragile, examples of
Sections~\ref{se:mve} and~\ref{se:boe} 
are not included in these results.

\subsection{Iteration Counts}\label{se:ic}
For each type of problem  we ran tests of the algorithms in
dimensions $3 \times 2^n$ for $n=0,1,2,3,4,5,6,7$.
Dimension 3 is our starting point since the seven algorithms
degenerate to three simple and effective algorithms in
dimension 2.
We use the factor 2 increase to sample higher dimensions
with less frequency as we get higher.  We believe this yields
a reasonable sample of low, intermediate and high dimensional 
problems.  

Note that a colourful feasibility problem instance in dimension $d$ consists of
$(d+1)^2$ points in dimension $d$.  Thus the size of
the input is cubic in $d$.  At present it is logistically difficult
to generate and store a colourful feasibility problem in dimension $d=1,000$.
After dimension 100, it also becomes increasingly difficult
to cope with numerical errors, especially for the algorithms
that include nearest point calculations, namely \B, \Bp\ and \H.
For this reason we do not include results for these algorithms
beyond $d=96$ for except for the relatively well-conditioned 
\rg\ problems where we stopped at $d=192$.

As one would expect, the guess-and-check algorithm \Rand\ 
performs badly as $d$ increases, except on problems from the
\Sp\ generator which have an abundance of solutions.  
We only include results from the \Rand\ algorithm when they
can be completed in a reasonable amount of time.

\ifapps
  The results of our computational experiments are presented
  in the graphs below and the tables in Appendix~\ref{ap:tables}.
\else
  The results of our computational experiments are presented
  in the graphs below.  We have made the tables containing the
  raw data for these graphs available at \cite{Advol}.
\fi
Each graph presents results for
a single random generator on a log-log scale with the average
iteration count of each algorithm plotted against the dimension.  
Additionally, the tables contain the values of the largest
iteration count observed in each type of trial; these show
the similar trends to the averages, although we notice that
\BO\ and \BOp\ sometimes perform substantially worse than the 
average, especially in the presence of ill-conditioning. 
The reasons for this are discussed in Section~\ref{se:boe}.

\begin{figure}[h!bt]
\begin{center}
\includegraphics[width=16cm,height=10.2cm]{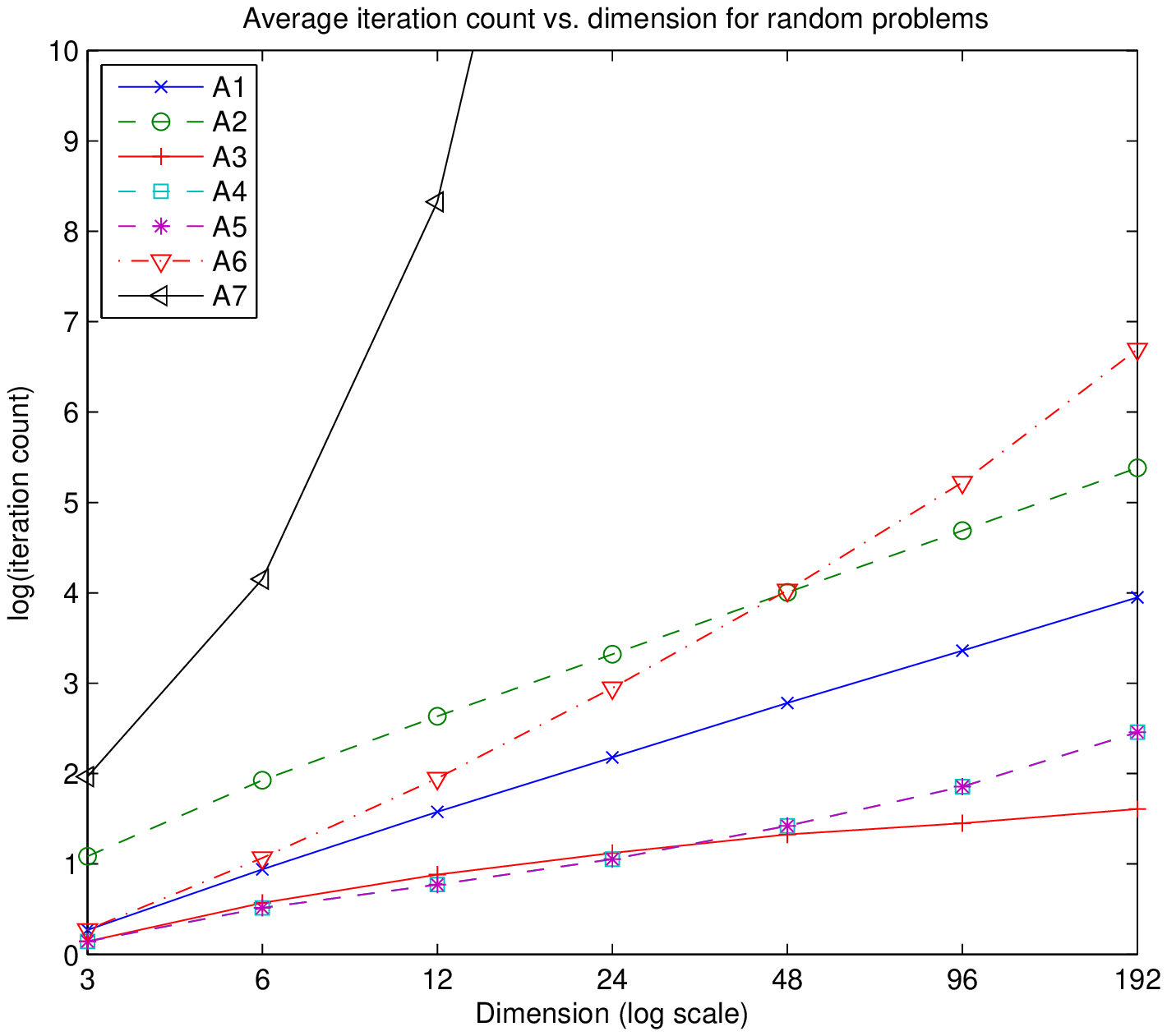}
\caption{Results for \rg.} \label{fig:random}
\end{center}
\end{figure}

\begin{figure}[h!bt]
\begin{center}
\includegraphics[width=16cm,height=10.2cm]{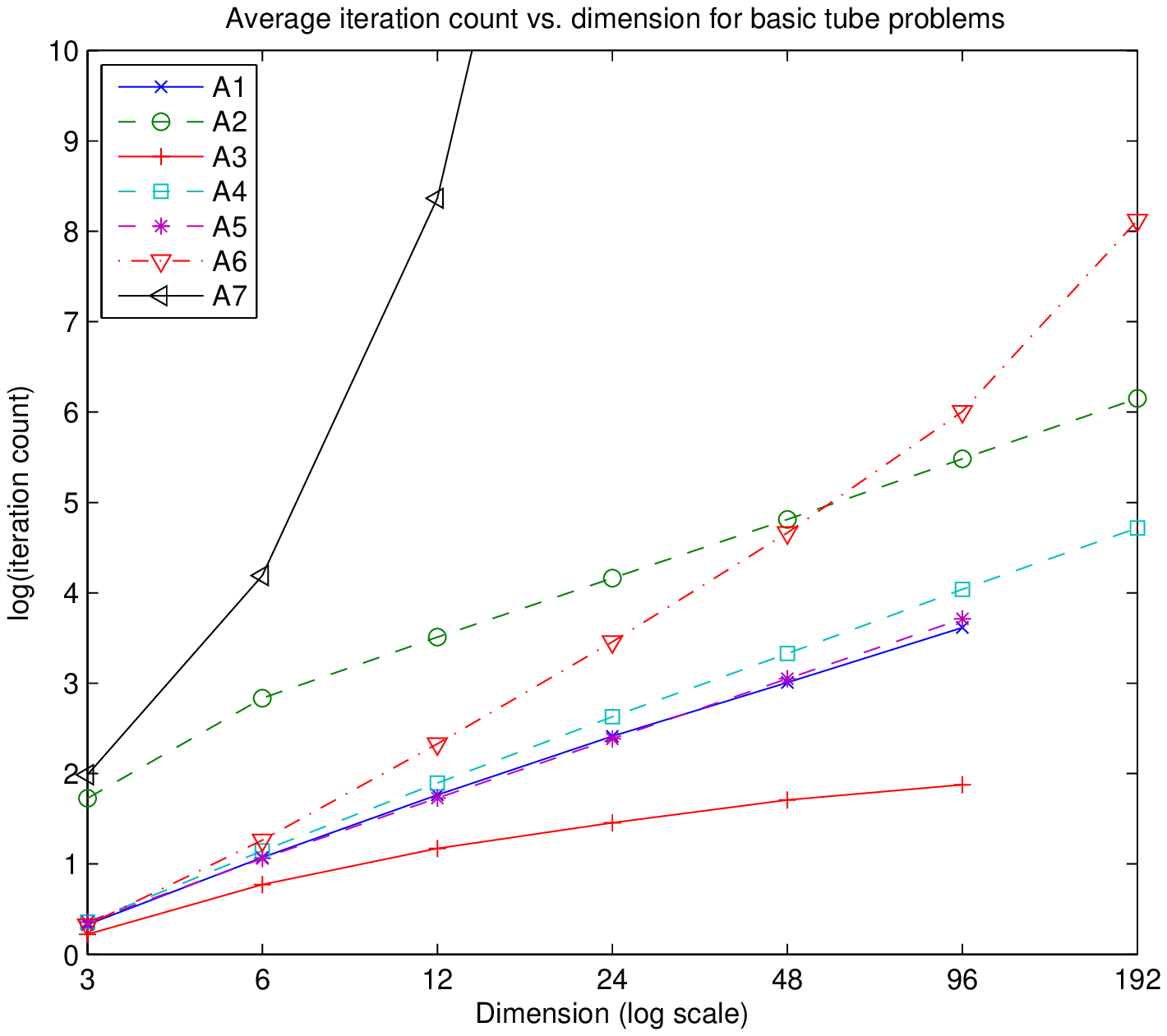}
\caption{Results for \tz.} \label{fig:tz}
\end{center}
\end{figure}

\begin{figure}[h!bt]
\begin{center}
\includegraphics[width=16cm,height=10.2cm]{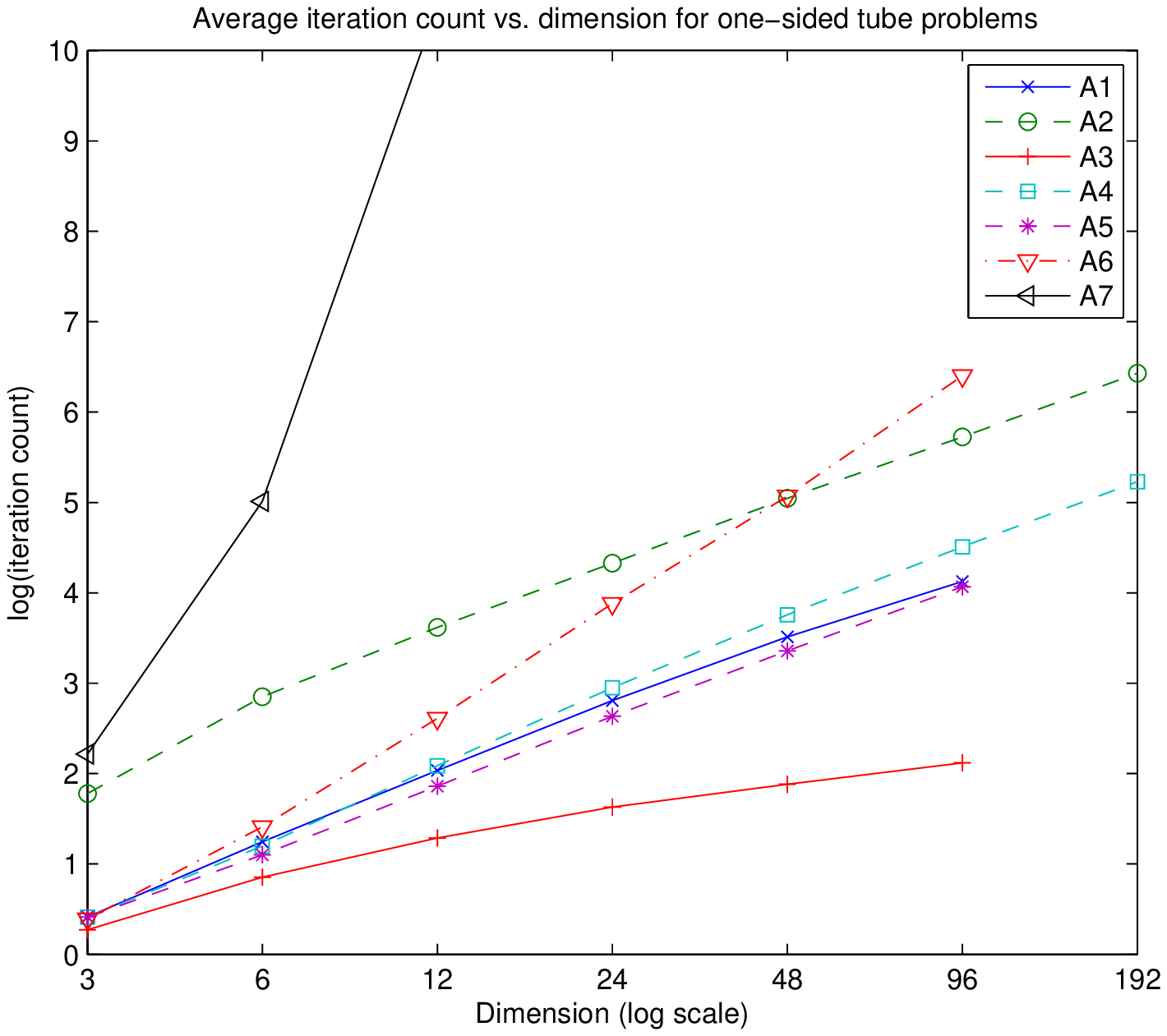}
\caption{Results for \to.} \label{fig:to}
\end{center}
\end{figure}

\begin{figure}[h!bt]
\begin{center}
\includegraphics[width=16cm,height=10.2cm]{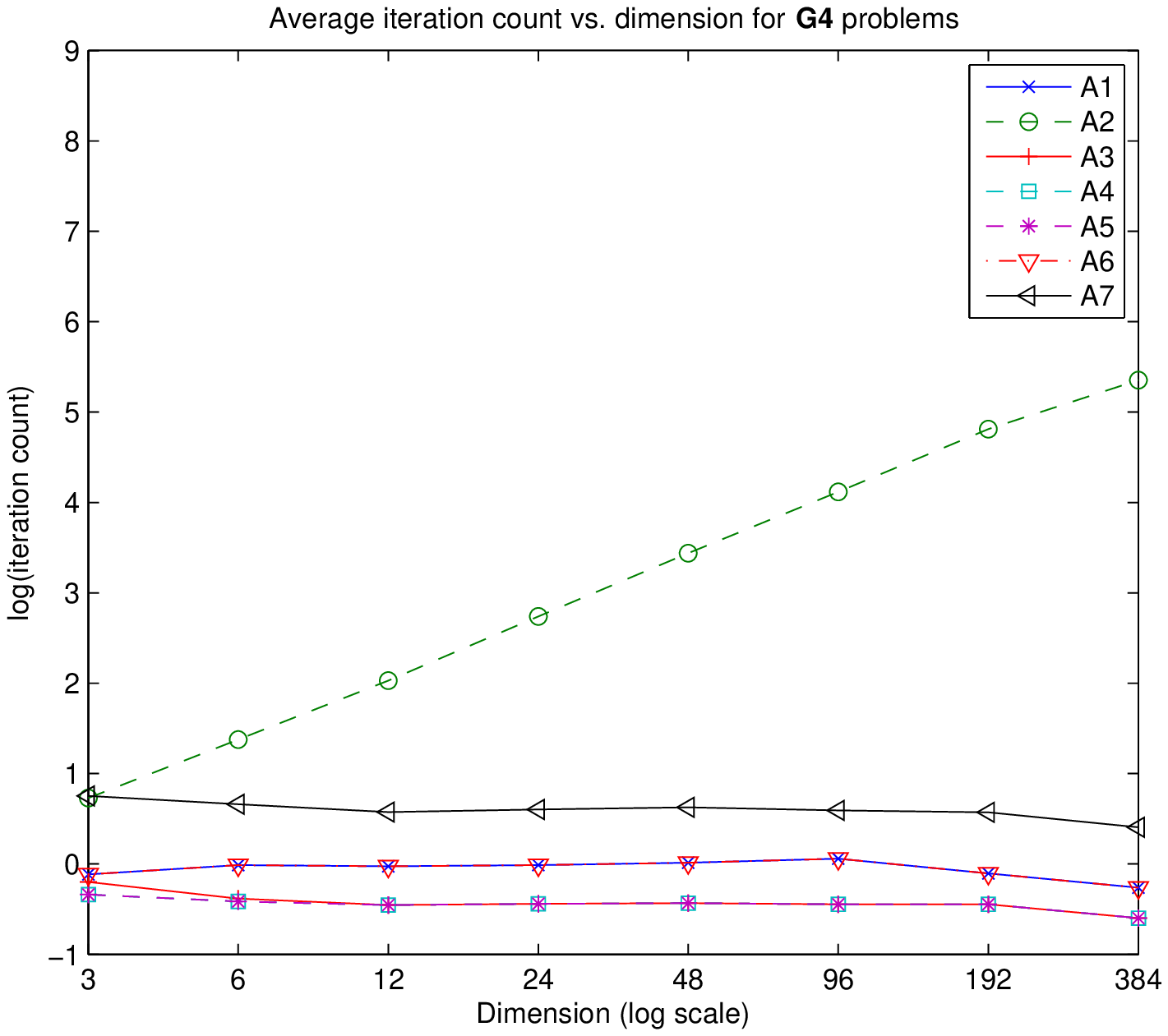}
\caption{Results for \Sp.} \label{fig:sp}
\end{center}
\end{figure}

\begin{figure}[h!bt]
\begin{center}
\includegraphics[width=16cm,height=10.2cm]{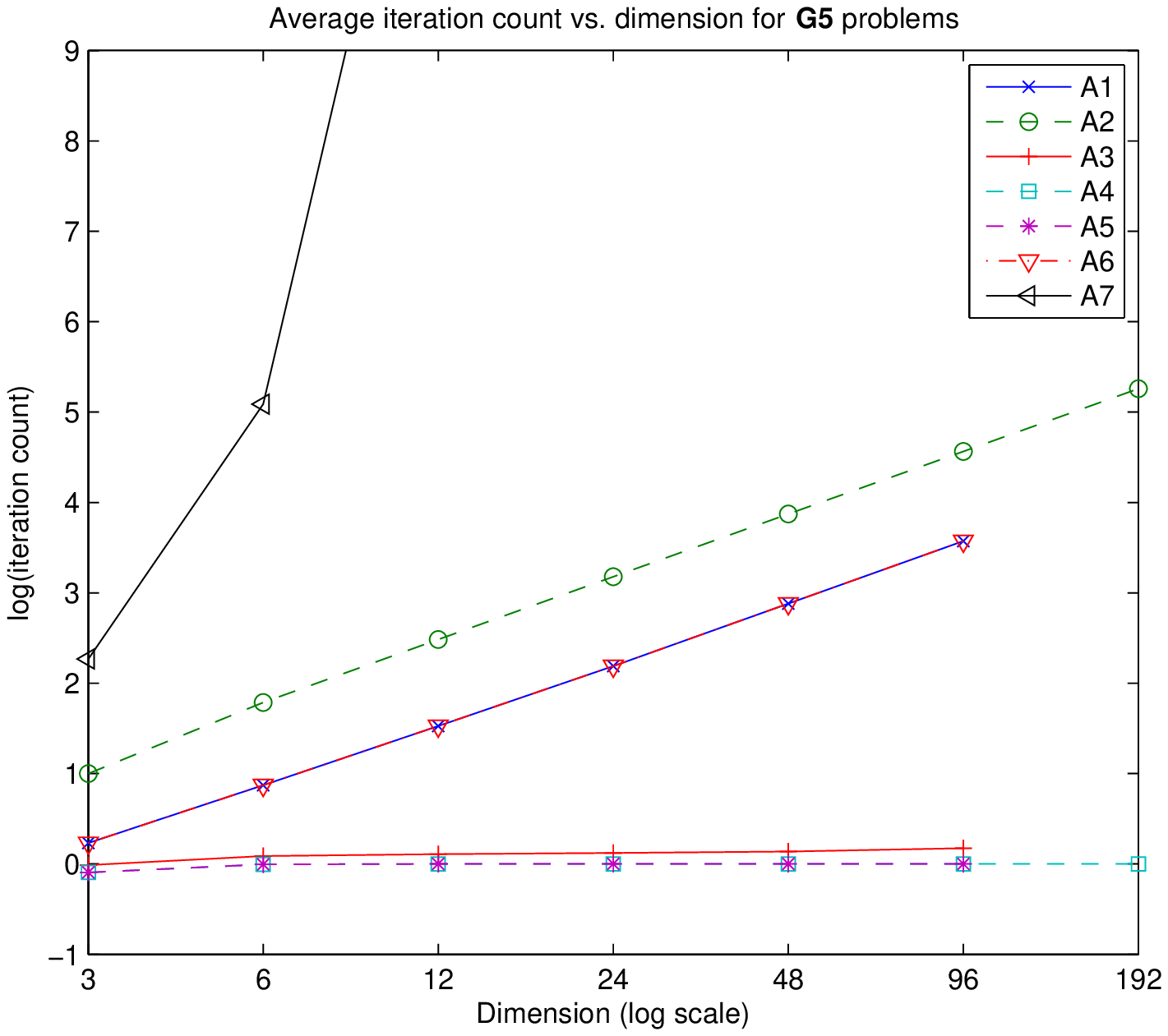}
\caption{Results for \Sd.} \label{fig:sd}
\end{center}
\end{figure}

\begin{figure}[h!bt]
\begin{center}
\includegraphics[width=16cm,height=10.2cm]{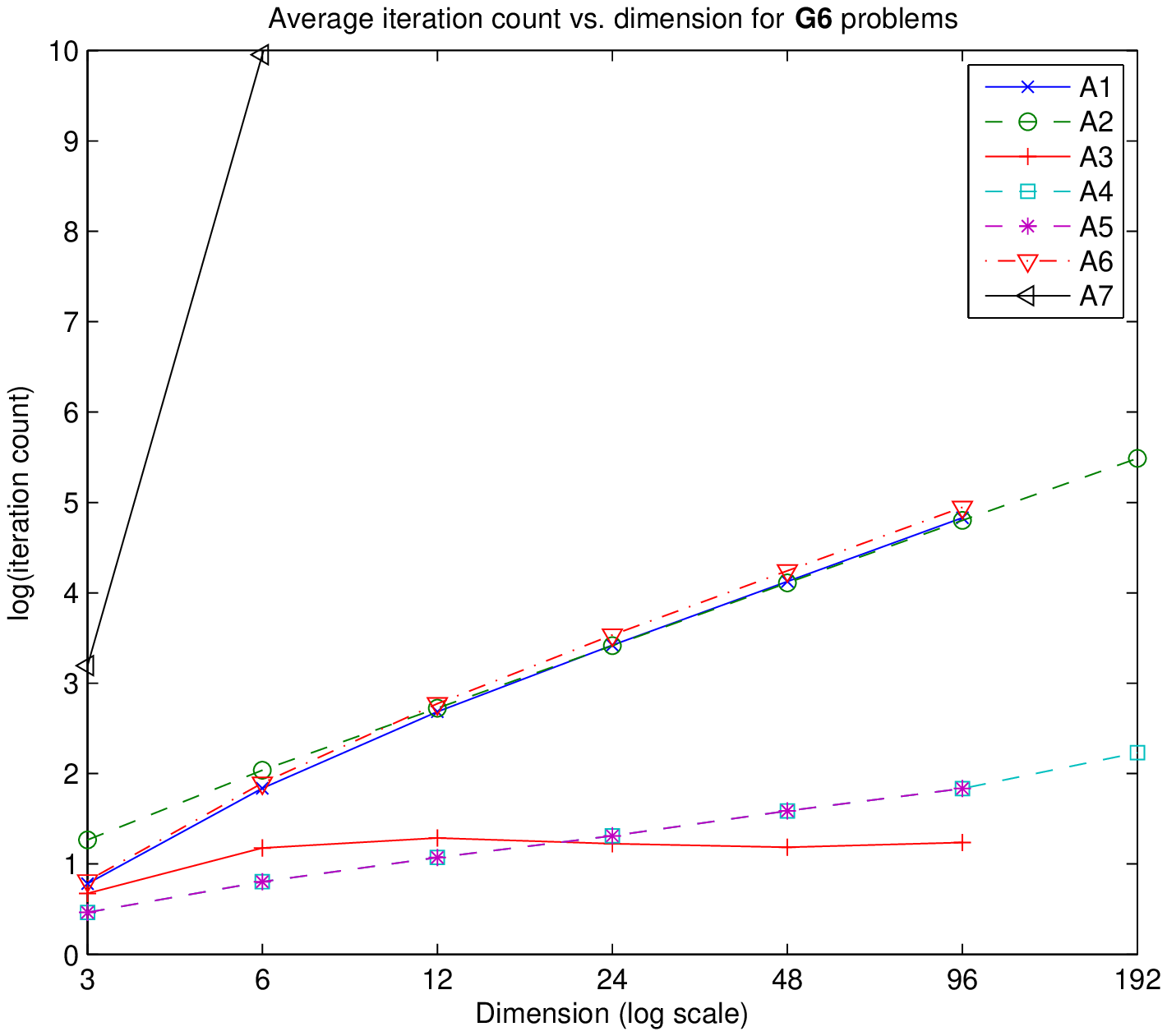}
\caption{Results for \Sm.} \label{fig:sm}
\end{center}
\end{figure}

For each generator at $d=3$ we sampled 100,000 problems,
at $d=6$ and $d=12$ we sampled 10,000 problems, at $d=24$ and
$d=48$ we sampled 1,000 problems and finally for $d \ge 96$
we sampled 100 problems.  
\ifapps
  Because of the varying sample sizes,
  it may not be entirely fair to compare the maxima listed in
  Appendix~\ref{ap:tables} between dimensions.
\fi
% All plots are generated using Sui Huang's MATLAB code \cite{Hua05}.
The results are plotted on as log-log graphs in 
Figures~\ref{fig:random}--\ref{fig:sm}.
We remark that polynomials appear asymptotically linear in
log-log plots, with the slope of the asymptote being the exponent
of the leading term of the polynomial and the $y$-intercept of
the asymptote representing the lead coefficient.

In Figure~\ref{fig:random} we see that \B\ and \BO\ appear to 
be taking a polynomial number of iterations to solution, while \MV\ and
\Rand\ do not appear to be polynomial.  Since each algorithm 
takes a polynomial time per iteration, the graphs of time versus
dimension show similar trends.

For the tube experiments, we used an angle parameter of
$\pi/6$, which is to say that all the vectors used made
an angle of at most $\pi/6$ with the $x$-axis.  
Smaller angles produce worse results for \BO, \BOp\ 
and \MV.  The example of \MV\ cycling, see Section~\ref{se:mve}
and Appendix~\ref{ap:mve},
was found using a smaller angle with \tz.

The tube experiments summarized in Figures~\ref{fig:tz} 
and~\ref{fig:to} show the impact of ill-conditioning on
all the algorithms.  For \B, \Bp, \H\ and \MV, convergence is
slightly slower and numerical errors become more common.
With these algorithms, our experiments began to
crash at dimension 192.  By contrast for the better conditioned
problems from \rg, the three algorithms with minimum distance
calculations crashed only at dimension 384 and \MV\ 
would in any case take too long on problems of this size.
Nevertheless, these algorithms remain effective at $d=96$.

The algorithms \BO\ and \BOp\ are more robust in the sense
that they are not as prone to crashes due to numerical errors.
This is the advantage of relying entirely on straightforward
linear algebra computations rather than considering nearest points
or volumes.
At the same time, they converge much more slowly due to problems 
of the type described in Section~\ref{se:boe} and Appendix~\ref{se:boe}.  

If we decrease the angle parameter which controls
the width of the tube and hence the conditioning, the results 
become more pronounced.  That is to say, \B, \Bp, \H\ and \MV\ 
become less stable numerically and experience a
further mild degradation in performance when not affected
by numerical errors, while \BO\ and \BOp\ become substantially
slower.

We comment that the \Rand\ algorithm performs about the same on \tz\ 
problems as it did on \rg\ problems.  This simply means that \tz\ 
problems typically have a similar number of solutions to \rg\ 
problems.  As one would expect, solutions to the one-sided tube problems
generated by \to\ are rarer than solutions to \rg\ and \tz\ problems
since the most of the points are clustered on one side.
Hence \Rand\ performs much worse on this type of
problem.

The problems with many solutions produced by \Sp\ 
are solved very quickly by all the algorithms, as
illustrated in Figure~\ref{fig:sp}.  In this case
the random sampling algorithm \Rand\ offers excellent
performance.  With the abundance of solutions, most of
the algorithms solve such problems in an expected 
constant number of iterations.  The exception is \BO\ 
which needs $\Theta(d)$ iterations at the
start to unwind the nearest point substitute $y$ from
a vertex to an interior point on a facet.
Since all the algorithms begin by checking the feasibility
of the initial simplex, the \Sp\ problems are often
solved in 0 iterations.

For the simplex structured problems of \Sd, we see all
the algorithms except \Rand\ perform very well, despite
the relative scarcity of solutions.  
We see that the other algorithms have
exactly the proper response to this structure -- they
systematically take points near vertices that are not
part of the current set.
In the case of \B, a new vertex of the simplex will be 
added at each step to give convergence in at most $d$
iterations, for \BO\ it takes one pass through the
$d+1$ colours, and for the multi-update algorithms 
\Bp, \BOp\ and \H\ one or two passes through the colours.
Algorithm \MV\ also solves these problems in a reasonable
number of iterations.

Finally, we see that the problems from \Sm\ 
where solutions are scarce are indeed more difficult
than random problems, but that, except for the \Rand\ 
algorithm, the impact on algorithmic performance is mild.
See Figure~\ref{fig:sm}.  
Curiously, the \Sm\ problems are the most difficult
problems for the \B\ algorithm.  % At the same time,
The multi-update algorithms \Bp, \BOp\ and \H\ 
perform extremely well.

\subsection{Cost per Iteration}\label{se:tpi}
\begin{figure}[h!bt]
\begin{center}
\includegraphics[width=16cm,height=10.2cm]{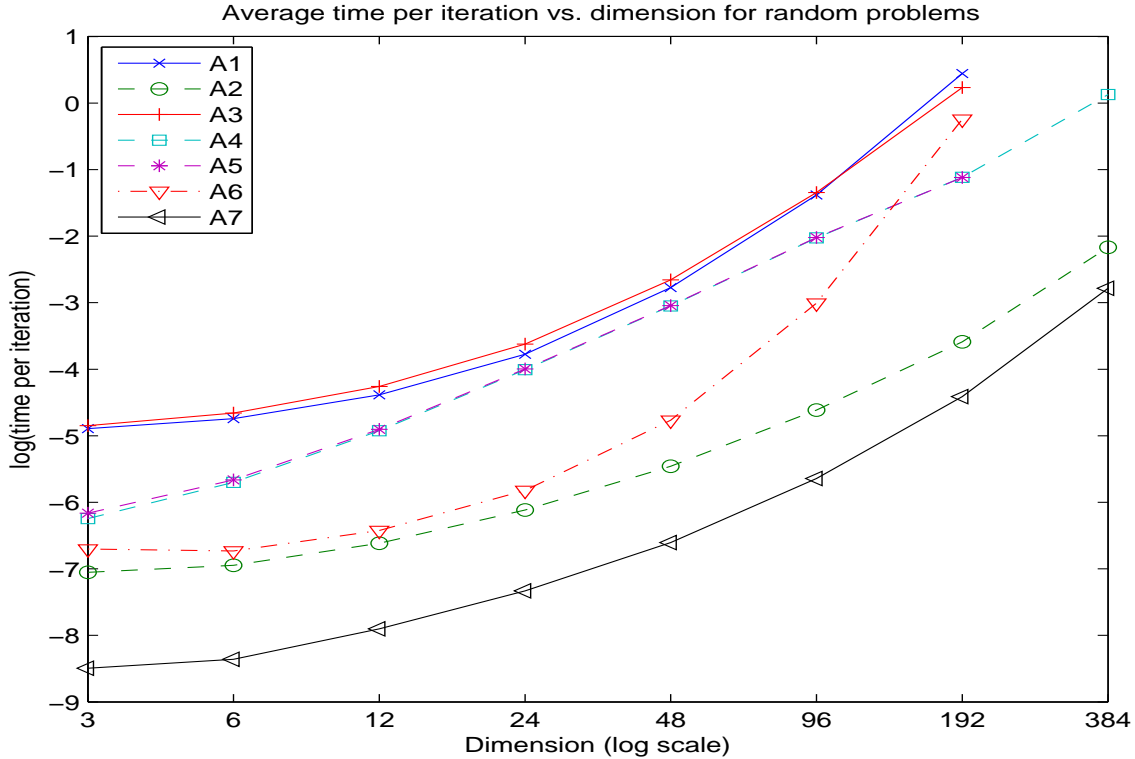}
\caption{Average iteration time of the algorithms.} \label{fig:tpi}
\end{center}
\end{figure}

In Figure~\ref{fig:tpi} we present
the average iteration times observed for all seven algorithms
on problems from the \rg\ generator.  
\ifapps
  The raw data for this graph is in Appendix~\ref{ap:tpi}.
\fi
We comment that the average time to complete an iteration does 
not change significantly with the problems type, so we have not
included the similar graphs for other generators.  
The data shows that in our implementation of these algorithms, 
the average time for an iteration is never very large.
For the slowest algorithms in the highest dimensions
the average iteration took less than 2 seconds.  

We see some interesting trends in the graphs.  
First, in low dimensions all the iteration times are
very fast and are presumably dominated by fixed startup
costs.  As the dimension increases, we begin to see
the asymptotic behaviour.  The algebraic algorithms \BO\ and \BOp\ 
show the expected $\Theta(d^3)$ behaviour, which appears linear
in the log-log plot.  
Asymptotically, the average time for an iteration of \BOp\ 
is about 10 times longer for an iteration of \BO.  

The algorithms \B\ and \Bp, which depend on a minimum distance
calculation, take longer on average to complete an iteration than \BOp.
The extra cost for the multiple updates in \Bp\ 
is relatively small.  However, the asymptotic slope of
these lines appear higher than for \BO, which means that
the nearest point calculations are causing the iterations
to take time $\Omega(d^3)$.  
The algorithm \MV\ has iteration times not much worse than \BO\ 
in low dimension, but its asymptotics look 
close to $O(d^4)$ as suggested in Section~\ref{se:mv}.
Algorithm \Rand\ exhibits $\Theta(d^3)$ iteration time and
is asymptotically about twice as fast on average per iteration
than \BO.

Unlike the other algorithms, the average iteration time for \H\ 
will be substantially affected by the conditioning of the problem.
Using the well-conditioned \rg\ problems, \H\ usually degenerates 
to \BOp\ and has a very similar average iteration time.
As the problems become more ill-conditioned, \H\ will begin
to use \Bp\ steps as well, and the average iteration time will 
increase towards the average iteration time for \Bp.

%
% CONCLUSIONS
%
\section{Conclusions}\label{se:conclusions}
Our experiments reveal several features of colourful
feasibility algorithms.  After considerable searching,
we found a problem instance which caused \MV\ to cycle.
We also found that \BO\ and \BOp\ can converge extremely
slowly in the face of ill-conditioning although \B\ and \Bp\
continue to perform reasonably well on the same 
examples.  We conclude that computationally the best algorithms are 
\Bp\ and \BOp\ and remark that these tightened algorithms
do yield substantial gains over the originals.

\subsection{A Cycling Example for \MV\ in Dimension 4}\label{se:mve}
% \ifapps % Now include this appendix in all versions.
  In Appendix~\ref{ap:mve} we exhibit an example in dimension 4
  for which the maximum volume heuristic cycles.
% \else
%   In \cite{Advol} we exhibit an example in dimension 4
%   for which the maximum volume heuristic cycles.
% \fi
This example was found using our tube generator \tz\ 
to produce configurations where for each colour, four points are
tightly bunched around (-1,0,0,0,0) and the fifth point
is close to (1,0,0,0,0) or vice-versa.  
The example is fairly ill-conditioned,
but not excessively so:  we rounded the values we found
for text formatting purposes, 
and observed that $\zero$ remained in the core and
that the behaviour of the algorithm was unaffected.

Close examination of the iterations of this example turns
up nothing out of the ordinary.  Since this example shows
that \MV\ can cycle, it is remarkable that it happens
so rarely.  It did not occur in the entire test suite of
Section~\ref{se:results}.  We tested extensively in dimensions
3 and 4, and were unable to find any examples of cycling in
dimension 3 or any examples of cycling in dimension 4 with
cycle length shorter than 6.  Higher dimensions and longer
cycle lengths do occur.  

One explanation for the results is that as one might expect, \MV\ 
is an effective heuristic in a typical situation.  
The distinguishing feature of the few bad examples is that
the points are placed in such a way that the simplices cluster
into a few groups of similar shape and volume.  
The heuristic of taking the maximum volume is then not very
helpful in choosing promising simplices.  
We note that this example is solved easily by the other algorithms.

\subsection{Flip-flopping During Convergence for \BO: 40,847 Iterations 
in Dimension 3}\label{se:boe}
We constructed an example of a colourful feasibility problem 
in dimension 3 that takes 40,847 iterations to solution
using a basic implementation of \BO.
% \ifapps  % Now this appendix is included in all versions.
  The exact points we used are contained in Appendix~\ref{ap:boe}.
% \else
%    The exact points we used are contained in \cite{Advol}.
% \fi
The algorithm is initialized with the simplex that uses
the first point of each colour.  At the
fifth iteration, the algorithm reaches a situation where the
current point $y$ lies on a facet $F$ of colours 2, 3 and 4 very close
to $\zero$.  Using this point the algorithm will pick the
point of colour 1 that has minimum dot product with $y$.
The second and third points of colour 1 lie almost in the directions of $y$
and $-y$, however neither of these forms a simplex with $F$
containing $\zero$.  In fact the fourth point of colour 1 does 
form a simplex containing $\zero$ with $F$, but it is nearly 
orthogonal to $y$.  As a result,
after two iterations, \BO\ returns to the same simplex.
The point $y$ will be recomputed at each step, and is slightly
closer to $\zero$ when the algorithm returns to the previous
simplex.  However, the improvement is quite small.  
Of course $\rho$ is also very small, so this is consistent with 
the performance guarantee described in Section~\ref{se:bo}.
The algorithm then proceeds to return to the same simplex
more than 20,000 times, with an incremental improvement to
$y$ at each iteration before finally taking the fourth
point of colour 1 and terminating.

As one would expect with a very ill-conditioned problem, this
example is numerically fragile -- the current version of our code
normalizes the coordinates before starting and does not suffer 
the same fate.  However bad behaviour is fairly typical. 
The tube generator for ill-conditioned problems in \cite{Hua05}
produces problems whose ill-conditioning depends on a parameter
defining the width of the tube.  As the width decreases, we
get an increasing number of cases where \BO\ and \BOp\ take
enormous numbers of iterations.

We remark that, in contrast, \B\ never returns to the same
simplex, so it cannot suffer from this
type of flip-flopping.  Indeed in dimension 3 it could do
no worse than visiting all $4^4=256$ simplices.  
At least 10 of these must contain $\zero$,
see \cite{BM06}, so the algorithm must terminate in at
most 246 iterations.  It is quite hard to see how this
limit could be approached.  The authors wonder if
a Klee-Minty-like example, see \cite{KM72}, 
of worst-case behaviour for \bara's pivoting
algorithm could be constructed.  

\subsection{Overall Effectiveness of Algorithms}\label{se:effective}
Despite the examples of Sections~\ref{se:mve} and~\ref{se:boe},
the results presented in Section~\ref{se:results} show that,
except for \Rand\ and to a lesser degree \MV, all the algorithms 
did a good job of solving all the problems.
We did find that the methods which include nearest point
calculations were more vulnerable to numerical errors than
\BO\ and \BOp, since our implementations began to crash
once we got much past $d=100$, especially on ill-conditioned
problems.  For the most part, reduced iteration counts of
the nearest point algorithms do not offset the extra time
spent per iteration compared to \BO\ and \BOp, since neither
iteration count is very high.  
In some cases of extreme ill-conditioning, such as in
Section~\ref{se:boe}, \BO\ and \BOp\ will take many
additional iterations and be much slower compared to the
nearest point algorithms.
In this situation either a hybrid algorithm such as \H,
or the basic \B\ or \Bp\ would work better.

We had hoped that the hybrid algorithm \H\ would offer
the benefits of \BOp, namely speed and robustness in
high dimensions, while stopping long periods of 
flip-flopping from occurring.  This did happen to a degree,
but in our benchmarking experiments the net time savings
were negligible, while \H\ retained \Bp's tendency to
crash due to numerical errors in high dimension.

\subsection{Advantages of Multiple Updates and Initialization}\label{se:multi}
The multi-update algorithms \Bp\ and \BOp\ do provide
substantial gains over their single update counterparts,
\B\ and \BO.  In the case of \Bp, we get a large reduction
in iteration count at very little cost in terms of iteration
time.  In our benchmarking experiments, this produced times
that were competitive with \BO\ and much better than \B.
The gains for \BOp\ relative to \BO\ are less impressive.
In our benchmarking experiments, \BOp\ consistently
averaged a 10\% to 40\% savings in total time to solution.

In fact, \BO\ is not as well suited as \B\ to take
advantage of multiple updates.  The point $y$ close to 
$\zero$ computed by \BO\ will almost always lie in the
interior of a facet of $\Deltab$, meaning that \BO\ 
will only have a single candidate colour to pivot.
In contrast, in high dimension, the closest point $x$ to 
$\zero$ will often lie on 
a relatively low dimensional face of $\Deltab$, allowing
multiple updates throughout the algorithm.

One difficulty for \BO\ is that it 
begins with $y$ at a vertex.  In a normal situation,
the first $d$ steps of \BO\ will each increase the
dimension of the smallest face containing $y$ by one
until $y$ lies in the interior of a facet, without
necessarily yielding a much better current simplex.
The multi-update \BOp\ does this all in the first iteration
in less time than it takes \BO\ to do $d$ steps.

We have not discussed the effects of the initial simplex
in this paper, but we can employ various heuristics to
choose a good initial simplex.  A few of these are implemented
in \cite{Hua05}.  We found that the most useful initialization 
heuristic was to run the first iteration of \BOp.  This runs
in $O(d^3)$ time and improves the subsequent iteration counts of
the algorithms, with the obvious exception of \Rand.
Even \BOp\ experiences a reduced iteration count, since the
point $y$ found by the initialization is not passed to the
algorithm.

\subsection{Theoretical Complexity of the Algorithms}\label{se:theory}
In Section~\ref{se:algs}, we remarked that \baraonn\ 
proved a worst-case bound for \B\ and \BO\ of $O(1/\rho^2)$
iterations up to numerical considerations and we improved
their iteration time for \BO\ from $O(d^4)$ to $O(d^3)$.  
We also mentioned that
we do not expect the multi-update and hybrid algorithms 
to improve the theoretical bounds.  
From the example of Section~\ref{se:mve}, we see that \MV\ 
is not guaranteed to converge.  The expected running time
of \Rand\ is 1 over the probability that random simplex
contains $\zero$, i.e.~around $2^d$ for random problems,
and as bad as $(d+1)^{d+1}/(d^2+1)$ for the type of
problems generated by \Sm.

The poor performance of \BO\ 
on ill-conditioned problems and examples like that
of Section~\ref{se:boe} confirm the worst-case
predictions of \baraonn's analysis.
On the other hand, we did not see this type of
behaviour for \B, and it is hard to see how
it could occur.

The model proposed in Section~\ref{se:sg} is that a
pure pivoting algorithm such as \B, defines a set of 
rooted trees on the $(d+1)^{d+1}$ simplices.
Each simplex which contains $\zero$ is the root of a tree, and
we draw an edge between the vertices representing simplices
$\Deltab_1$ and $\Deltab_2$ if when \B\ encounters $\Deltab_1$
it pivots to $\Deltab_2$.  Then the worst performance of the
algorithm in terms of the number of iterations would be the height
of the highest tree.  A smart algorithm will produce short
trees by pivoting several simplices to a given simplex at
a lower level.  

Consider a situation where trees have a constant expansion
factor $k$ near the base, that is, low level vertices are
connected to roughly $k$ vertices in the level above.
The number of trees is $p(d+1)^{d+1}$ where $p$ is the 
probability that a simplex contains $\zero$.  
If the trees expand up to height $h$, each
tree will contain on the order of $k^h$ vertices.
Then we must have $k^h p(d+1)^{d+1} \le (d+1)^{d+1}$,
the total number of vertices.  Rearranging, we get
$h \le -\log_k(p)$.  This expression predicts the average
iteration count for \B\ to grow linearly for \rg\ 
problems, to be constant for \Sp\ problems and to
grow at $\Theta(d \log d)$ for \Sm\ problems.
All of these match very well with our observed results.
The \Sd\ problems are predicted to be more difficult than
they are observed to be, but that is not surprising given
their simple structure.

\subsection{Future Considerations}\label{se:future}
We finish by returning to the motivating question of \baraonn:
Is there a polynomial time algorithm for colourful feasibility?
By improving the implementation of \BO, we have improved the
worst case for this algorithm from $O(d^4/\rho^2)$ to
$O(d^3/\rho^2)$, however the dependence on $\rho$ has
not improved.  Indeed our experiments give strong evidence
that the analysis for \BO\ is tight.

The situation for \B\ is less clear.  We do not see the
same bad behaviour with ill-conditioned problems that
we found for \BO, so it is possible that a better
guarantee exists for this algorithm.  
In light of the model suggested in Section~\ref{se:theory}
it is quite difficult to see how to construct a
Klee-Minty-like bad case for \B\ as discussed in 
Section~\ref{se:boe}.  We view this as an appealing challenge.

%
% ACKNOWLEDGMENTS
%
\section{Acknowledgments}
This research was supported
by NSERC Discovery grants for the four authors,
by the Canada Research Chair program for the first and last authors
and by a MITACS grant for the second and third authors.
The third author worked on this project as part of the
Discrete Optimization project of the IMO at the University of Magdeburg.

%
% APPENDICES
%
\appendix

%
% APPENDIX A: Max_vol example.
%
\section{Example in dimension 4 where \MV\ cycles}\label{ap:mve}
This example consists of 5 normalized points in each of the
5 colours in $\R^4$.
The points are presented in Table~\ref{ta:mvc}.
They are grouped by colour, with the rows
representing $x$, $y$, $z$ and $w$ coordinates, respectively.

\begin{table}[h!bt]
%\begin{tabular}{rc}
 Red points 

 \begin{tabular}{|c|c|c|c|c|} \hline
  -0.98126587 &  0.99234170 & -0.99375618 & -0.98428021 & -0.99649986 \\ \hline
  0.13481464 &  0.01125213 & -0.01676635 & -0.03542019 &  0.03152825 \\ \hline
  0.00569666 & -0.12300509 &  0.10203928 &  0.17121850 &  0.07625092 \\ \hline
  0.13751313 &  0.00104048 & -0.04189897 &  0.02494182 & -0.01340880 \\ \hline
 \end{tabular} \vspace{0.5mm} \\
 
 Green points 

 \begin{tabular}{|c|c|c|c|c|} \hline
  0.99924734 & -0.99225276 &  0.95301586 &  0.99770745 &  0.98808067 \\ \hline
  0.03530047 & -0.07048563 &  0.17760263 &  0.03405179 & -0.00874509 \\ \hline
  -0.01500068 &  0.10036231 & -0.24516979 & -0.01526145 & -0.12973853 \\ \hline
  0.00579663 & -0.01984027 &  0.01048096 &  0.05645716 &  0.08238952 \\ \hline
 \end{tabular} \vspace{0.5mm} \\
 
 Blue points 

 \begin{tabular}{|c|c|c|c|c|} \hline
  -0.98758195 & -0.99742900 & -0.97286388 & -0.97433105 &  0.99536963 \\ \hline
  -0.03897365 &  0.02836725 &  0.13575382 &  0.14413058 & -0.06519965 \\ \hline
  -0.14957699 & -0.06348104 & -0.17638005 &  0.17286629 &  0.06380946 \\ \hline
  -0.02810110 & -0.01734511 & -0.06322067 & -0.00475659 &  0.03027639 \\ \hline
 \end{tabular} \vspace{0.5mm} \\
 
 Tan points  

 \begin{tabular}{|c|c|c|c|c|} \hline
  0.99782436 &  0.99917562 &  0.95584087 & -0.98768930 &  0.96962649 \\ \hline
  0.01692290 &  0.03972232 &  0.17806542 & -0.10337937 &  0.14481818 \\ \hline
  0.03437294 & -0.00816965 & -0.21878711 &  0.09313650 & -0.12491250 \\ \hline
  0.05365310 &  0.00186470 &  0.08242045 & -0.07147128 &  0.15247636 \\ \hline
 \end{tabular} \vspace{0.5mm} \\
 
 White points 

 \begin{tabular}{|c|c|c|c|c|} \hline
  -0.99979855 & -0.97268376 & -0.97231627 & -0.95622769 &  0.99791825 \\ \hline
  0.00600345 &  0.06950105 &  0.21172943 & -0.29221243 & -0.02997771 \\ \hline
  0.00415788 & -0.00409898 & -0.03733932 & -0.01550644 &  0.01616939 \\ \hline
  0.01869548 &  0.22144776 &  0.09152860 &  0.00022801 & -0.05476362 \\ \hline
 \end{tabular}
%\end{tabular}
\caption{Coordinates of points of an example where \MV\ cycles in dimension 4.}
\label{ta:mvc}
\end{table}

The initial simplex is taken to be (1,1,1,1,1), i.e.,~the first point
of each colour.  The algorithm proceeds to visit simplices
(1,1,4,1,1), (3,1,4,1,1), (3,1,4,3,1), (3,1,1,3,1) and (1,1,1,3,1)
before returning to the original simplex and repeating.

%
% APPENDIX B: Flip-flopping example.
%
\section{Example in dimension 3 where \BO\ takes 40,847 iterations}
\label{ap:boe}
This example consists of 4 unnormalized points in each of the
4 colours in $\R^3$.
The points are presented in Table~\ref{ta:bol}.
They are grouped by colour, with the rows
representing $x$, $y$ and $z$ coordinates, respectively.
\begin{table}[h!bt]
Red points 

\begin{tabular}{|c|c|c|c|} \hline
1.00000320775369 & -0.01000436049274 & -0.01000129525998 &  1.00000089660284 \\
\hline
0.00000340785030 &  0.99999739350954 & -1.00000497855619 &  0.00000051797159 \\
\hline
0.00999859615603 &  0.00000371775824 &  0.00000030149139 & -0.01999639732055 \\
\hline
\end{tabular}

Green points 

\begin{tabular}{|c|c|c|c|} \hline
1.00000363763560 & -0.00999644886160 & -0.00999943004295 &  1.00000335962280 \\
\hline
-0.00000325123594 &  1.00000064545156 & -1.00000169806216 & -0.00000080450760 \\
\hline
0.01000493174811 & -0.00000024088601 &  0.00000009099437 & -0.01999811804365 \\
\hline
\end{tabular}

Blue points 

\begin{tabular}{|c|c|c|c|} \hline
0.99999949817337 & -0.00999587145461 & -0.00999627213896 &  0.99999551963712 \\
\hline
-0.00000260397964 &  1.00000485455718 & -1.00000419710665 & -0.00000024626161 \\
\hline
0.00999854691703 &  0.00000123671997 & -0.00000381812529 & -0.01999801526314 \\
\hline
\end{tabular}

Tan points 

\begin{tabular}{|c|c|c|c|} \hline
0.99999980645233 &  0.10000000280522 & -0.60000327600988 &  0.99999642880542 \\
\hline
0.00000024487465 & -0.98999719313413 &  0.79999695643245 & -0.00000429109491 \\
\hline
0.01000455311709 & -0.00000405877812 &  0.00000372117690 & -0.01000272055280 \\
\hline
\end{tabular}
\caption{Coordinates of points of an example taking 40,847 iterations of \BO\ 
in dimension 3.}
\label{ta:bol}
\end{table}
 
The initial simplex is taken to be (1,1,1,1), i.e.,~the first point
of each colour.  It then updates to (1,3,1,1), (1,3,2,1), (1,3,2,3),
(1,3,2,2) and reaches (3,3,2,2) on the fifth iteration.  
At this point, it begins to flip
between (3,3,2,2) and (2,3,2,2) with $y$ initially alternating
between values close to (0.2,$\pm$0.00200,0.00285).
The values of all these coordinates decrease very slowly as the
algorithm continues.  At iteration 40,847 it chooses fourth point
of colour 1 instead of the third.  This makes the current simplex
(4,3,2,2) which contains $\zero$.

%
% APPENDIX C: Raw data.
%
\section{Iteration counts from our experiments}\label{ap:tables}
In this Appendix we present the raw data from our computational
experiments. % graphed in Section~\ref{se:ic}.
Each table presents results for
a single random generator.  The entries give the average number 
of iterations to solution for each algorithm at the given dimension.
For each generator at $d=3$ we sampled 100,000 problems,
at $d=6$ and $d=12$ we sampled 10,000 problems, at $d=24$ and
$d=48$ we sampled 1,000 problems and finally for $d \ge 96$
we sampled 100 problems.

\begin{table}[h!bt]\label{ta:rga}
\begin{tabular}{|l||c|c|c|c|c|c|c|} \hline 
~ & \B & \BO & \Bp & \BOp & \H & \MV & \Rand \\ \hline \hline
$d=3$  & 1.31 & 2.96 & 1.15 & 1.15 &  1.15 &  1.31 &      7.15 \\ \hline
$d=6$  & 2.56 & 6.87 & 1.77 & 1.67 &  1.67 &  2.90 &      63.48 \\ \hline
$d=12$  & 4.84 & 13.93 & 2.42 & 2.16 & 2.16 &  7.01 &     4133.15 \\ \hline
$d=24$  & 8.84 & 27.70 &  3.07 & 2.87 & 2.87 & 19.07 &    Large \\ \hline
$d=48$  & 16.14 & 54.88 & 3.77 & 4.14 & 4.14 & 56.12 &    Large \\ \hline
$d=96$  & 28.80 & 108.71 & 4.26 & 6.39 & 6.39 & 185.57 &  Large \\ \hline
$d=192$  & 51.96 & 217.59 & 4.99 & 11.68 & 11.68 & 808.78 & Large \\ \hline
$d=384$  & Unstable & 425.26 & Unstable & 21.63 & Unstable & Large & Large 
\\ \hline
\end{tabular}
\caption{Average iteration counts in \rg\ generator tests.}
\end{table}

\begin{table}[h!bt]\label{ta:rgm}
\begin{tabular}{|l||c|c|c|c|c|c|c|} \hline 
~ & \B & \BO & \Bp & \BOp & \H & \MV & \Rand \\ \hline \hline
$d=3$  & 5  & 136  &   4   & 4  & 4  & 5  &  102 \\ \hline
$d=6$  & 7  & 21  &    5   & 5  & 5  & 12  & 579 \\ \hline
$d=12$  & 10 & 30 &     6  &  6 &  6 &  20 &  47362 \\ \hline
$d=24$  & 15 & 37 &     6  &  8 &  8 &  43 &  Large \\ \hline
$d=48$  & 22 & 67 &     6  &  9 &  9 &  105 & Large \\ \hline
$d=96$  & 39 & 120 &    6  &  10 & 10 & 269 & Large \\ \hline
$d=192$  & 63 & 241 &    7  &  19 & 19 & 1574 & Large \\ \hline
$d=384$  & Unstable & 472 & Unstable & 30 & Unstable & Large & Large 
\\ \hline
\end{tabular}
\caption{Maximum iteration counts found in \rg\ generator tests.}
\end{table}

\begin{table}[h!bt]\label{ta:tza}
\begin{tabular}{|l||c|c|c|c|c|c|c|} \hline 
~ & \B & \BO & \Bp & \BOp & \H & \MV & \Rand \\ \hline \hline
$d=3$  & 1.39  & 5.62 &    1.25  & 1.43  &   1.43  & 1.38  &  7.30  \\ \hline
$d=6$  & 2.92  & 17.00 &   2.17  & 3.14  &   2.89  & 3.54  &  66.02 \\ \hline
$d=12$  & 5.83  & 33.48 &   3.23  & 6.65  &   5.64  & 10.26  & 4296.66 \\ \hline
$d=24$  & 11.18 & 64.30 &   4.29 &  13.86 &   10.86 & 31.75 &  Large \\ \hline
$d=48$  & 20.24 & 123.02 &  5.51 &  27.91 &   21.11 & 106.11 & Large \\ \hline
$d=96$  & 37.12 & 240.49 &  6.54 &  56.70 &   40.91 & 406.10 & Large \\ \hline
$d=192$  & Unstable & 468.52 & Unstable & 111.84 & Unstable & 3367.60 & Large 
\\ \hline
$d=384$  & Unstable & 909.82 & Unstable & 220.50 & Unstable & Large & Large 
\\ \hline
\end{tabular}
\caption{Average iteration counts in \tz\ generator tests.}
\end{table}

\begin{table}[h!bt]\label{ta:tzm}
\begin{tabular}{|l||c|c|c|c|c|c|c|} \hline 
~ & \B & \BO & \Bp & \BOp & \H & \MV & \Rand \\ \hline \hline
$d=3$  & 5  & 4783  & 4  & 5  & 5  & 6  & 109 \\ \hline
$d=6$  & 8  & 2880  & 6  & 44  & 10  & 14  & 1079 \\ \hline
$d=12$  & 13 & 842 & 8 & 60 & 14 & 33 & 78418 \\ \hline
$d=24$  & 21 & 217 & 9 & 36 & 23 & 78 & Large \\ \hline
$d=48$  & 31 & 249 & 9 & 55 & 41 & 258 & Large \\ \hline
$d=96$  & 47 & 323 & 9 & 77 & 76 & 840 & Large \\ \hline
$d=192$  & Unstable & 561 & Unstable & 140 & Unstable & 11784 & Large 
\\ \hline
$d=384$  & Unstable & 1013 & Unstable & 260 & Unstable & Large & Large 
\\ \hline
\end{tabular}
\caption{Maximum iteration counts found in \tz\ generator tests.}
\end{table}

\begin{table}[h!bt]\label{ta:toa}
\begin{tabular}{|l||c|c|c|c|c|c|c|} \hline 
~ & \B & \BO & \Bp & \BOp & \H & \MV & \Rand \\ \hline \hline
$d=3$  & 1.51 &  5.93   & 1.31 & 1.51   & 1.51 & 1.48  &  9.16 \\ \hline
$d=6$  & 3.48 &  17.26  & 2.35 & 3.31   & 3.01 & 4.10  &  150.31 \\ \hline
$d=12$ & 7.64 &  37.22  & 3.62 & 8.06   & 6.43 & 13.61 &  Large \\ \hline
$d=24$ & 16.59 & 75.73  & 5.11 & 19.11  & 13.92 & 48.51 & Large  \\ \hline
$d=48$ & 33.51 & 155.48 & 6.57 & 42.81  & 28.70 & 159.29 & Large \\ \hline
$d=96$ & 61.97 & 306.64 & 8.32 & 90.98  & 58.44 & 602.07 & Large \\ \hline
$d=192$ & Unstable &  619.55 & Unstable & 186.86 & Unstable & Large & Large
\\ \hline
$d=384$  & Unstable & 1221.43 & Unstable & 382.10 & Unstable & Large & Large 
\\ \hline
\end{tabular}
\caption{Average iteration counts in \to\ generator tests.}
\end{table}

\begin{table}[h!bt]\label{ta:tom}
\begin{tabular}{|l||c|c|c|c|c|c|c|} \hline 
~ & \B & \BO & \Bp & \BOp & \H & \MV & \Rand \\ \hline \hline
$d=3$ & 6 & 2756 & 5  & 6  &   6  & 6  &   127 \\ \hline
$d=6$ & 9 & 3704 & 7  & 38  &  9  & 14  &  1709 \\ \hline
$d=12$ & 16 &  689 & 8 & 55 &   16 & 46 &   Large \\ \hline
$d=24$ & 28 &  195 & 9 & 52 &   27 & 124 &  Large \\ \hline
$d=48$ & 50 &  257 & 10 & 83 &  47 & 505 &  Large \\ \hline
$d=96$ & 78 &  374 & 11 & 133 & 83 & 2023 & Large \\ \hline
$d=192$ & Unstable & 736 & Unstable & 226 & Unstable & Large & Large 
\\ \hline
$d=384$  & Unstable & 1399 & Unstable & 454 & Unstable & Large & Large 
\\ \hline
\end{tabular}
\caption{Maximum iteration counts found in \to\ generator tests.}
\end{table}

\begin{table}[h!bt]\label{ta:spa}
\begin{tabular}{|l||c|c|c|c|c|c|c|} \hline 
~ & \B & \BO & \Bp & \BOp & \H & \MV & \Rand \\ \hline \hline
$d=3$  & 0.89 & 2.07  & 0.82 & 0.71 & 0.71 & 0.89 & 2.12 \\ \hline
$d=6$  & 0.99 & 3.96  & 0.68 & 0.66 & 0.66 & 0.99 & 1.94 \\ \hline
$d=12$ & 0.97 & 7.61  & 0.63 & 0.63 & 0.63 & 0.97 & 1.78 \\ \hline
$d=24$ & 0.99 & 15.46  & 0.64 & 0.64 & 0.64 & 0.99 & 1.83 \\ \hline
$d=48$ & 1.01 & 31.15  & 0.65 & 0.65 & 0.65 & 1.01 & 1.87 \\ \hline
$d=96$ & 1.06 & 61.44  & 0.64 & 0.64 & 0.64 & 1.06 & 1.81 \\ \hline
$d=192$ & 0.90 & 122.88 & 0.64 & 0.64 & 0.64 & 0.90 & 1.77 \\ \hline
$d=384$  & 0.77 & 211.20 & 0.55 & 0.55 & 0.55 & 0.77 & 1.50
\\ \hline
\end{tabular}
\caption{Average iteration counts in \Sp\ generator tests.}
\end{table}

\begin{table}[h!bt]\label{ta:spm}
\begin{tabular}{|l||c|c|c|c|c|c|c|} \hline 
~ & \B & \BO & \Bp & \BOp & \H & \MV & \Rand \\ \hline \hline
$d=3$ & 2 & 5 & 2 & 3 & 3 & 2 &  38 \\ \hline
$d=6$ & 3 & 7 & 2 & 2 & 2 & 3 &  17 \\ \hline
$d=12$ & 6 & 12 & 1 & 1 & 1 & 6 & 30 \\ \hline
$d=24$ & 6 & 24 & 1 & 1 & 1 & 6 & 19 \\ \hline
$d=48$ & 5 & 48 & 1 & 1 & 1 & 5 & 16 \\ \hline
$d=96$ & 5 & 96 & 1 & 1 & 1 & 5 & 14 \\ \hline
$d=192$ & 3 & 192 & 1 & 1 & 1 &  4 & 15 \\ \hline
$d=384$  & 4 & 384 & 1 & 1 & 1 & 4 & 9
\\ \hline
\end{tabular}
\caption{Maximum iteration counts found in \Sp\ generator tests.}
\end{table}

\begin{table}[h!bt]\label{ta:sda}
\begin{tabular}{|l||c|c|c|c|c|c|c|} \hline
~ & \B & \BO & \Bp & \BOp & \H & \MV & \Rand \\ \hline \hline
$d=3$  & 1.26 & 2.72 &   0.99 & 0.91 & 0.91 & 1.26 &  9.67 \\ \hline
$d=6$  & 2.39 & 5.97 &   1.09 & 0.99 & 0.99 & 2.39 &  161.93 \\ \hline
$d=12$ & 4.61 & 12.00 &  1.12 & 1.00 & 1.00 & 4.61 &  Large \\ \hline
$d=24$ & 8.94 & 24.00 &  1.13 & 1.00 & 1.00 & 8.94 &  Large \\ \hline
$d=48$ & 17.82 & 48.00 & 1.15 & 1.00 & 1.00 & 17.82 & Large \\ \hline
$d=96$ & 35.58 & 96.00 & 1.19 & 1.00 & 1.00 & 35.58 & Large \\ \hline
$d=192$ & 71.15 & 192.00 & 1.47 & 1.00 & 1.00 & 71.15 & Large
\\ \hline
\end{tabular}
\caption{Average iteration counts in \Sd\ generator tests.}
\end{table}

\begin{table}[h!bt]\label{ta:sdm}
\begin{tabular}{|l||c|c|c|c|c|c|c|} \hline
~ & \B & \BO & \Bp & \BOp & \H & \MV & \Rand \\ \hline \hline
$d=3$ & 3  & 5   & 3 & 2 &  2 & 3 & 128 \\ \hline
$d=6$ & 5  & 6   & 3 & 1 &  1 & 5 & 1371 \\ \hline
$d=12$ & 9  & 12   & 3 & 1 &  1 & 9 & Large \\ \hline
$d=24$ & 14 & 24  & 2 & 1 & 1 & 14 & Large \\ \hline
$d=48$ & 24 & 48  & 2 & 1 & 1 & 24 & Large \\ \hline
$d=96$ & 41 & 96  & 2 & 1 & 1 & 41 & Large \\ \hline
$d=192$ & 81 & 192 & 3 & 1 & 1 &   81 & Large
\\ \hline
\end{tabular}
\caption{Maximum iteration counts found in \Sd\ generator tests.}
\end{table}

\begin{table}[h!bt]\label{ta:sma}
\begin{tabular}{|l||c|c|c|c|c|c|c|} \hline
~ & \B & \BO & \Bp & \BOp & \H & \MV & \Rand \\ \hline \hline
$d=3$  & 2.19 & 3.54 &   1.96 & 1.59 & 1.59  & 2.26 & 24.39 \\ \hline
$d=6$  & 6.27 & 7.67 &   3.24 & 2.23 & 2.23  & 6.65 & 21041.05 \\ \hline
$d=12$ & 14.64 & 15.23 & 3.63 & 2.92 & 2.92 & 16.03 & Large \\ \hline
$d=24$ & 30.55 & 30.42 & 3.40 & 3.71 & 3.71 & 34.25 & Large \\ \hline
$d=48$ & 61.96 & 60.95 & 3.27 & 4.89 & 4.89 & 69.65 & Large \\ \hline
$d=96$ & 125.31 & 121.73 & 3.45 & 6.26 & 6.26 & 140.79 & Large \\ \hline
$d=192$ & Unstable & 242.06 & Unstable & 9.31 & Unstable & Unstable & Large
\\ \hline
\end{tabular}
\caption{Average iteration counts in \Sm\ generator tests.}
\end{table}

\begin{table}[h!bt]\label{ta:smm}
\begin{tabular}{|l||c|c|c|c|c|c|c|} \hline
~ & \B & \BO & \Bp & \BOp & \H & \MV & \Rand \\ \hline \hline
$d=3$ & 5  & 7 & 5   & 4  & 4  & 6  & 242 \\ \hline
$d=6$ & 12 & 15 & 7  & 6 &  6 &  12 & 173941 \\ \hline
$d=12$ & 25 & 25 & 8  & 9 &  9 &  25 & Large \\ \hline
$d=24$ & 47 & 49 & 9  & 13 & 13 & 51 & Large \\ \hline
$d=48$ & 101 & 94 & 13  & 22 & 22 & 95 & Large \\ \hline
$d=96$ & 154 & 174 & 6  & 35 & 35 & 183 & Large \\ \hline
$d=192$ & Unstable & 331 & Unstable & 69 & Unstable & Unstable & Large
\\ \hline
\end{tabular}
\caption{Maximum iteration counts found in \Sm\ generator tests.}
\end{table}

\clearpage

%
% APPENDIX D: Times.
%
\section{Average time per iteration}\label{ap:tpi}
In Table~\ref{ta:tpi} we give the average CPU time per iteration
for our \rg\ experiments.  This was computed using the MATLAB
{\tt cputime} function.  

\begin{table}[h!bt]
\begin{tabular}{|l||c|c|c|c|c|c|c|} \hline
~ & \B & \BO & \Bp & \BOp & \H & \MV & \Rand \\ \hline \hline
$d=3$ & 0.0075 & 0.0009 & 0.0078 & 0.0019 & 0.0021 & 0.0012 & 0.0002 \\ \hline
$d=6$ & 0.0087 & 0.0010 & 0.0095 & 0.0033 & 0.0035 & 0.0012 & 0.0002 \\ \hline
$d=12$ & 0.0124 & 0.0013 & 0.0141 & 0.0073 & 0.0074 & 0.0016 & 0.0004 \\ \hline
$d=24$ & 0.0229 & 0.0022 & 0.0267 & 0.0182 & 0.0184 & 0.0030 & 0.0007 \\ \hline
$d=48$ & 0.0625 & 0.0043 & 0.0702 & 0.0474 & 0.0477 & 0.0085 & 0.0014 \\ \hline
$d=96$ & 0.2510 & 0.0099 & 0.2608 & 0.1318 & 0.1324 & 0.0495 & 0.0035 \\ \hline
$d=192$ & 1.5592 & 0.0277 & 1.2623 & 0.3275 & 0.3268 & 0.7843 & 0.0121 \\ \hline
$d=384$  & Unstable & 0.1144 & Unstable & 1.1381 & Unstable & Unstable & 0.0619
\\ \hline
\end{tabular}
\caption{Average iteration times on \rg\ generator tests.}\label{ta:tpi}
\end{table}

The time per iteration is fairly constant across problem types
so we do not include data from the other generators.  One
difference that will occur is that \H\ will have a higher average 
iteration time as that \BOp\ for ill-conditioned problems.
In random problems, we rarely see slow convergence of \BOp\ 
so it is unnecessary to use the slower steps from \Bp.
With ill-conditioned problems the \Bp\ steps become more frequent
and increase the average time per iteration.

%
% BIBLIOGRAPHY
%
% Setup for bibtex
\bibliographystyle{hacked}
% End setup for bibtex.

\bibliography{refs}

\providecommand{\bysame}{\leavevmode\hbox to3em{\hrulefill}\thinspace}
\providecommand{\MR}{\relax\ifhmode\unskip\space\fi MR }
% \MRhref is called by the amsart/book/proc definition of \MR.
\providecommand{\MRhref}[2]{%
  \href{http://www.ams.org/mathscinet-getitem?mr=#1}{#2}
}
\providecommand{\href}[2]{#2}
\begin{thebibliography}{DHST05}

\bibitem[B{\'a}r82]{Bar82}
I.~B{\'a}r{\'a}ny, \emph{A generalization of {C}arath\'eodory's theorem},
  Discrete Math. \textbf{40} (1982), no.~2-3, 141--152.

\bibitem[BM05]{BM06}
I.~B{\'a}r{\'a}ny and J.~Matou{\v s}ek, \emph{Quadratically many colorful
  simplices}, submitted, 2005.

\bibitem[BO97]{BO97}
I.~B{\'a}r{\'a}ny and S.~Onn, \emph{Colourful linear programming and its
  relatives}, Math. Oper. Res. \textbf{22} (1997), no.~3, 550--567.

\bibitem[DHST05]{DHST06}
A.~Deza, S.~Huang, T.~Stephen, and T.~Terlaky, \emph{Colourful simplicial
  depth}, Discrete Comput. Geom. (2005), ~To appear.  
  ~{\tt arXiv:math.CO/0506003}

\bibitem[Hua]{Hua05}
S.~Huang, \emph{{MATLAB} code for colourful linear programming}, available
  at:\\ http://optlab.mcmaster.ca/{\textasciitilde}huangs3/CLP/ and \\
  http://www.math.uni-magdeburg.de/{\textasciitilde}stephen/Software/CLP/.

\bibitem[KM72]{KM72}
V.~Klee and G.~J. Minty, \emph{How good is the simplex algorithm?},
  Inequalities III, Proc. 3rd Symp., Los Angeles 1969, Academic Press, 1972,
  pp.~159--175.

\bibitem[ST05]{ST06}
T.~Stephen and H.~Thomas, \emph{A quadratic lower bound for colourful
  simplicial depth}, in preparation, 2005.

\bibitem[WW01]{WW01}
U.~Wagner and E.~Welzl, \emph{A continuous analogue of the upper bound
  theorem}, Discrete Comput. Geom. \textbf{26} (2001), no.~2, 205--219.

\end{thebibliography}

\end{document}